\documentclass[]{article}  %我们以后用这个就行
\pdfoutput=1
\usepackage{graphicx}
\usepackage{float}
\usepackage{amsthm}
\usepackage{amssymb}
\usepackage{amsmath}
\usepackage{mathrsfs}
\usepackage{indentfirst}
\usepackage{geometry}
\usepackage{authblk}
\newtheorem{Theorem}{Theorem}[section]
\theoremstyle{definition}

\newtheorem{Lemma}[Theorem]{Lemma}

\theoremstyle{definition}

\title{Near automorphisms of $G_{(n,m)}$}
\author{Songnian Xu, Wenhao Zhen, Dein Wong\thanks{Corresponding author. E-mail address: wongdein@163.com.}}
\affil{\textit{School of Mathematics, China University of Mining and Tecnology, Xuzhou, China.}}
\date{}

\begin{document}
\baselineskip 17pt

\title{Near automorphisms of $G_{(n,m)}$}

\author{Songnian Xu\\
{\small  Department of Mathematics, China University of Mining and Technology}\\
{\small Xuzhou, 221116, P.R. China}\\
{\small E-mail: xsn1318191@cumt.edu.cn}\\ \\
Wang\thanks{Corresponding author}\\
{\small Department of Mathematics, China University of Mining and Technology}\\
{\small Xuzhou 221116, P.R. China}\\
{\small E-mail: }}

\date{}
\maketitle

\begin{abstract}
Let $G$ be a graph with vertex set $V(G)$, $f$ a permutation of $V(G)$.
Define $\delta_f(G)=|d(x,y)-d(f(x),f(y))|$ and $\delta_f(G)=\Sigma\delta_f(x,y)$, where the sum is taken over all unordered pair $x$, $y$ of distinct vertices of $G$.
$\delta_f(x,U)=\Sigma\delta_f(x,y)$, where $U\subseteq V(G)$ and $y\in U$.
Let $\pi(G)$ denote the smallest positive value of $\delta_f(G)$ among all permutations of $V(G)$.
A permutation $f$ with $\delta_f(G)=\pi(G)$ is called a near automorphisms of $G$\cite{HV}.
In this paper, we define $G_{(n,m)}$ is a graph obtained from $K_n$ by add $t_i$ pendent vertices to $y_i$ which is a vertex of $K_n$, $i=1,\cdots,m$, and we say $y_i$ is a c-pendent vertex of $G_(n,m)$.
We determine $\pi(G_{(n,m)})$ and describe permutations $f$ of $G_{(n,m)}$ for which $\pi(G_{(n,m)})=\delta_f(G_{(n,m)})$.
Because $G_{(1,1)}$ is a star and it is easy, hence we let $n\geq 2$.
Suppose $G_(n,m)$ has $m$ c-pendent vertices $\{y_1, \ldots, y_m\}$ and $y_i$ has $t_i$ pendent vertices($1\leq t_1\leq t_2\leq \ldots \leq t_m$).
For $m<n$ we have
$$\pi(G_{(n,m)})=
\left\{
    \begin{array}{lc}
        2n-4 & n \leq t_1+2, m=1 \\
        2t_1&otherwise\\
    \end{array}
\right.
$$
\\
For $m=n$ we have
$$\pi(G_{(n,n)})=
\left\{
    \begin{array}{lc}
        4 & t_1=1,t_2=2 \\
        2t_1+2t_2&otherwise\\
    \end{array}
\right.
$$
\end{abstract}

\let\thefootnoteorig\thefootnote
\renewcommand{\thefootnote}{\empty}

\section{Introduction}
Let $G$ be a connected graph with vertex set $V(G)$.
The distance $d_G(x,y)$ between two vertices $x$ and $y$ is the length of a shortest $x-y$ path (if the graph is clear from the text, $d_G(x,y)$ will be abbreviated to $d(x,y)$).
For a permutation $f$ of $V(G)$, set

$\delta_f(x,y)=|d(x,y)-d(f(x),f(y))|$, $\delta_f(x)=\Sigma_{y\in V(G)}\delta_f(x,y)$; $\delta_f(x,U)=\Sigma_{y\in U}\delta_f(x,y)$ for $U\subseteq V(G)$.

The total relative displacement of $f$ is defined as $\delta_f(G)=\Sigma\delta_f(x,y)$, where the sum is taken over all unordered pairs $x$, $y$ of distinct vertices of $G$.
Clearly, $\Sigma_{x\in V(G)}\delta_f(x)=2\delta_f(G)$, and $f$ is an automorphism if and only if $\delta_f(G)=0$.
Let $\pi(G)$ denote the smallest positive value of $\delta_f(G)$ among all permutations $f$ of $V(G)$.
A $near$ $automorphism$ of $G$ is a permutation $f$ of $V(G)$ such that $\sigma_f(G)=\pi(G)>0$.
It is clear that a connected graph has a near automorphism if and only if $G$ it is not complete.
The displacement graph of $G$ with respect to $f$ is the directed multigraph $G[f]$ whose vertex set $V(G[f])=\{a_1,\ldots,a_t\}$, where $t=diam(G)$, and arc set $A(G[f])=\{\langle a_i,a_j\rangle: i\neq j$, there is a pair of vertices $u$ and $v$ such that $d(u,v)=i$ and $d(f(u),f(v))=j\}$.
And let $\alpha(u,v)=\langle a_i,a_j\rangle$, where $d(u,v)=i$ and $d(f(u),f(v))=j$.

Let $P_n$ (resp., $C_n$) be the path (resp., the cycle) with $n$ vertices. Chartrand et al. \cite{HV} conjectured that $\pi(P_n)=2n-4$.
Later , Aitken \cite{AW} verified this conjecture and characterized those permutation $f$ with $\delta_f(P_n)=2n-4$.
In 2008, Chang et al. \cite{CCF} proved that $\pi(C_n)=4\lfloor n/2\rfloor-4$ and characterized the near automorphism of $C_n$.
Another interesting result along this line was given by Reid \cite{R} which determined $\pi(K_{n_1,n_2,\ldots,n_r})$.
Chang and Fu \cite{CF} gave characterizations for trees $T$ with $\pi(G)=2$ and trees $T$ with $\pi(T)=4$.
Till now, not much is known for chaotic mappings.

In this paper, we defined $G_{(n,m)}$ is a graph obtained from $K_n$ by add $t_i$ pendent vertices to $y_i$ which is a vertex of $K_n$, $i=1,\cdots,m$, and we say $y_i$ is a c-pendent vertex of $G_(n,m)$. And we proved $\pi(G_{n,m})=2t_1$, $2t_1+2t_2$, $4$, or $2n-4$ and characterized the near automorphisms $G_(n,m)$.

\section{Preliminaries}
\begin{Lemma}\cite[Lemma 3]{CF}
If $G$ is a graph and $f$ is a permutation of $V(G)$ which is not a automorphism, then there is an edge $(u,v)$ of $G$ such that $\delta_f(u,v)\geq1$
\end{Lemma}

\begin{Lemma}\cite[Lemma 6]{CF}
For each vertex $a_i\in V(G[f])$, $1\leq i\leq$ diam$(G)$, $deg^{+}(a_i)=deg^{-}(a_i)$.
\end{Lemma}

\begin{Lemma}\cite[Lemma 3.2]{R}
If $N=p_1+p_2+\cdots+p_k$ is a partition of integer $N\geq2$ into $k\geq2$ positive integers $p_1,p_2,\ldots,p_k$, then $\Sigma p_rp_s\geq(3(k-1)/(k+1))(N-1)$, where the sum is over all $k \choose 2$ pair chosen from ${p_1,p_2,\ldots,p_k}$.
Moreover, equality occurs if and only either $k=2$ and $p_1=1$ and $p_2=N-1$ (or $p_1=N-1$ and $p_2=1$) or $k=N=3$ (so that $p_1=p_2=p_3=1$).
\end{Lemma}
\section{near automorphisms of $G_{(n,m)}$}
In this article, we are going to assume $U=U_1\cup U_2 \cup\cdots\cup U_m$ and $W=\{y_1,y_2,\ldots,y_m\}$ which is the set of c-pendent vertices of $G_{(n,m)}$.
\begin{Lemma} We let $f\neq1$ is a permutation of $G_{(n,m)}$ and $n> m\geq2$.
Then $\delta_f(G)>2t_1$ if $f(U)\neq U$.
\end{Lemma}
\begin{proof}
Assume $x_{ij}\in U_i$ and $f(x_{ij})\in W$, then we have $d(f(x_{ij}),f(x_{ks}))=1$ or $2$ for $i\neq k$ and $d(x_{ij},x_{ks})=3$, hence $\delta_f(x_{ij})\geq t_1(3-2)=t_1$.
It is easy to prove that $\delta_f(G)\geq 2t_1+2>2t_1$ by Lemma 2.1 and Lemma 2.2.
\end{proof}

\begin{Lemma} We let $f\neq1$ is a permutation of $G_{(n,m)}$, $n> m\geq2$ and $2\leq|U_1|=t_1\leq|U_2|=t_2\leq\cdots\leq|U_m|=t_m$.
If $f(U_i)=f(M_{i_1})\cup f(M_{i_2})\cup\cdots \cup f(M_{i_h})$ and $h\geq2$ then $\sigma_f(G)>2t_1$.
\end{Lemma}

\begin{proof}
$Case$ $1$. If there is another $U_j$ such that $f(U_j)=f(M_{j_1})\cup f(M_{j_2})\cup\cdots \cup f(M_{j_{h_1}})$ $(h_1>2)$ where $f(M_{j_r})\subseteq U_{j_r}$ and $\{j_1,j_2,\ldots,j_{h_1}\}\subseteq\{1,2,\ldots,m\}$.
 We prove $\delta_f(G)\geq2t_i+2t_j$.

Since $d(x_{ij_1},x_{ij_2})=2$ and $d(x_{ij},x_{ks})=3$ for $j_1\neq j_2$ and $i\neq k$, $\delta_f(U_i)=\Sigma_{i_j\neq i_k}|M_{i_j}||M_{i_k}|$.

By Lemma 2.3, $\delta_f(U_i)=\Sigma_{i_j\neq i_k}|M_{i_j}||M_{i_k}|\geq(3(h-1)/(h+1))(t_i-1)$.
Since $f(h)=(3(h-1)/(h+1))(t_i-1)$ is a monotonically increasing function of $h$, $\delta_f(U_i)$ take the minimum $t_i-1$ value if and only if $h=2$ and $|M_{i_1}|=t_i-1$ and $|M_{i_2}|=1$ (or $|M_{i_1}|=1$ and $|M_{i_2}|=t_i-1$) of $h=t_i=3$ (so that $|M_{i_1}|=|M_{i_2}|=|M_{i_3}|=1$) by Lemma 2.3.
On the other hand,  we can also get $\delta_f(U_j)$ take the minimum value $t_j-1$ if and only if $h_1=2$ and $|M_{j_1}|=t_j-1$ and $|M_{j_2}|=1$ (or $|M_{j_1}|=1$ and $|M_{j_2}|=t_j-1$) of $h_1=t_j=3$ (so that $|M_{j_1}|=|M_{j_2}|=|M_{j_3}|=1$).
Firstly, we consider the case where $h=h_1=2$.

Assume w.l.o.g. that
\begin{center}
$M_{i_1}=\{x_{i1}\}$, $M_{i_2}=\{x_{i2},\ldots,x_{it_i}\}$
\end{center}

\begin{center}
 $M_{j_1}=\{x_{j1}\}$, $M_{j_2}=\{x_{j2},\ldots,x_{jt_j}\}$
\end{center}
Since $\alpha(x_{i1},x_{i2})=\cdots=\alpha(x_{i1},x_{it_i})=\alpha(x_{j1},x_{j2})=\cdots=\alpha(x_{j1},x_{jt_j})=(a_2,a_3)$, $\alpha(x_{i_1},y_i)=(a_1,a_2)$ or $\alpha(x_{i_2},y_i)=(a_1,a_2)$ and $\alpha(x_{j_1},y_j)=(a_1,a_2)$ or $\alpha(x_{j_2},y_j)=(a_1,a_2)$, $\delta_f(G)\geq 2[(t_i-1)+(t_j-1)+2]=2t_i+2t_j$.
By observation and discussion above, we know that $\delta_f(G_{(n,m)})=2t_i+2t_j$ if and only if there is exactly one pair of vertices $(x_i,x_j)$ such that $\alpha(x_i,y_j)=\alpha(x_j,y_i)=(a_2,a_1)$ and $\delta_f(x,y_k)=0$ for $x\in U\backslash \{x_i, y_i\}$ and $k=1,2,\cdots,m$.

Now, we consider the case where $h=3$ or $h_1=3$.
We are going to use a same discussion to $h=h_1=2$, and we are going to get $\delta_f(G)\geq 2[(t_i-1)+(t_j-1)+3]=2t_i+2t_j+2>2t_i+2t_j$.
And this result will be useful in the discussion of $G_{(n,n)}$.

$Case$ $2$. $f(U_j)\subseteq U_{j'}$ for $j\neq i$ and $j, j'\in\{1,2,\ldots,m\}$.

If $h>2$ or $h=2$ and $|M_{i_1}|\neq1$ (or$|M_{i_2}|\neq1$), it is easy to find that $\delta_f(G)>2t_i$.
So we assume that $h=2$, $f(U_i)=f(M_{i_1})\cup f(M_{i_2})$ and $M_{i_1}=\{x_{i1}\}$, $M_{i_2}=\{x_{i2},\ldots,x_{it_i}\}$.
Since $\alpha(x_{i1},x_{i2})=\cdots=\alpha(x_{i1},x_{it_i})=(a_2,a_3)$, $\alpha(x_{i1},y_i)=(a_1,a_2)$ or $\alpha(x_{i2},y_i)=\alpha(x_{i3},y_i)=\cdots=\alpha(x_{it_i},y_i)=(a_1,a_2)$, $\delta_f(G)\geq2[\delta_f(U_i)+\delta_f(y_i,U_i)]\geq2t_i$.
If $t_i>t_1$, then we have $\delta_f(G)>2t_1$.
If $t_i=t_1$, then $|U_{i_2}|\geq|U_i|$.

If $|U_{i_2}|=|U_i|$.
So there is exactly one vertex $x_l\in U_l$ ($l\neq i$) where $f(x_l)\in U_{i_2}$.
But, it is obvious that $f(U_l)\nsubseteq U_{i_2}$, a contradiction.

If $|U_{i_2}|>|U_i|$.
So, there are at least two points $x,y\in U\backslash U_i$ such that $f(x),f(y)\in U_{i_2}$.
Since $\alpha(x,x_{i2})=\cdots=\alpha(x,x_{it_i})=\alpha(y,x_{i2})=\cdots=\alpha(y,x_{it_i})=(a_3,a_2)$,  $\delta_f(G)\geq2[(t_i-1)+(t_i-1)]+2=2t_i+2(t_i-1)>2t_i>2t_1$.

Let $M=U_{k+1}\cup U_{k+2}\cup\cdots\cup U_{m}$ and $T=U_1\cup U_2\cup\cdots\cup U_k$.
If $|U_i|\geq 2$ and $f(U_i)\neq U_i$, we know that there is one pair of vertices $(x_1,x_2)$ such that $\alpha(x_1,x_2)=(a_2,a_3)$ for $x_1,x_2\in M$.
And we have $\delta_f(G)>2$ by Lemma 2.1 and Lemma 2.2.
Next we assume $f(U_i)=U_i$ for $|U_i|\geq 2$.
Then there is a $\{x_i1\}=U_i$ such that $f(x_i1)=x_j1$ and $f(y_i)$.

If $M=U_{k+1}$, the result is obvious.
Now we assume that $m>k+1$, we know that $\sigma_f(G_{(n,m)})\geq 2t_{k+1}+2t_{t+2}>2$ if $f(M)=M$ by Case 1, so we have $f(x)\in T$ for some $x\in W$ and that also means that $f(y)\in W$ for some $y\in T$.
So all of this tells us that there is a $\{x_{j1}\}=U_j$ such that $f(x_{j1})\in U_l\subseteq M$.

Next we prove that $\delta_f(G)>2$.

Since $U_l\geq 2$, there is a $f(x)\in U_l$ and $\alpha(x_{j1},x)=(a_3,a_2)$.
Hence we have $\delta_f(G)>2$ by Lemma 2.1 and Lemma 2.2.
\\

$Case$ $3$. $1=|U_1|=|U_2|=\cdots=|U_m|$.

Since we need $\delta_f(G)>0$, there is a $\{x_{i1}\}=U_i$ such that $f(x_{i1})=x_{l1}$, $f(y_i)\neq y_l$ and $f(y)=y_l$ for $i\neq l$.
If $y\in K_n\setminus \{y_1,y_2,\ldots,y_m\}$, then there is a $y_q\in\{y_1,y_2,\ldots,y_m\}$ such that $f(y_q)\in  K_n\setminus \{y_1,y_2,\ldots,y_m\}$, so $\delta_f(G)\geq \delta_f(x_{i1},y_i)+\delta_f(x_{i1},y)+\delta_f(x_{q1},y_q)=1+1+1>2$.
If $y\in \{y_1,y_2,\ldots,y_m\}$, we also have $\delta_f(G)>2$.

To sum up, the Lemma is proved.

\end{proof}
\begin{Lemma} We let $f$ is a permutation of $G_{(n,m)}$ and $n>m\geq2$, then $\sigma_f(G)>2t_1$ if $f(y_i)=y_j$ for $i\neq j\in\{1,2,\ldots,m\}$.
\end{Lemma}

\begin{proof}. By Lemma 3.1, we assume that $f(U_i)=U_i$.
So $\delta_f(G)\geq\delta_f(y_i,U_i)+\delta_f(y_i,U_j)+\delta_f(y_j,U_j)=2t_j+t_i>2t_1$.
\end{proof}

\begin{Theorem} We let $f$ is a permutation of $G_{(n,m)}$ and $n> m\geq2$, then $\pi(G_{n,m})=2t_1$ and $\delta_f(G_{(n,m)})=2t_1$ if and only if $f(y)=y_i$, $f(U_l)=U_l$, $f(y_j)=y_j$ and $f(y_i)\in V(K_n\setminus \{y_1,y_2,\ldots,y_m\})$ where $y\in V(K_n\setminus \{y_1,y_2,\ldots,y_m\})$, $|U_i|=t_1$ and $i\neq j\in\{1,2,\ldots,m\}$.
\end{Theorem}

\begin{proof}. By Lemma 3.1 and Lemma 3.2, it is easy to know that $f(U_l)=U_l$ and $f(y_i)\neq y_j$ for $i\neq j$ and $l,i,j\in\{1,2,\ldots,m\}$.
And there must be some $y\in V(K_n\setminus \{y_1,y_2,\ldots,y_m\})$ that satisfies $f(y)=y_i$ and $f(y_i)\in V(K_n\setminus \{y_1,y_2,\ldots,y_m\})$, otherwise $\delta_f(G)=0$ or $\delta_f(G)>2t_1$.
Since $\delta_f(G)\geq\delta_f(y,U_i)+\delta_f(y_i,U_i)=2t_i\geq2t_1$, equality occurs if and only if $t_i=t_1$, $f(y_j)=y_j$ for $i\neq j\in\{1,2,\ldots,m\}$.
\end{proof}

Next we consider near automorphisms for $G_{(n,1)}$ (we also let $n\geq2$).

\begin{Theorem} We let $y_1$ be the only divergent vertex of $K_n$ and $U$ be the set of pendent vertices of $y_1$.
Then
$$\pi(G)=
\left\{
    \begin{array}{lc}
        2n-4 & n \leq |U|+2=t+2 \\
        2t&otherwise\\
    \end{array}
\right.
$$

\end{Theorem}

\begin{proof}. Assume $U=\{x_1,x_2,\ldots,x_t\}$ and $V(K_n)=\{y_1,y_2,\ldots,y_n\}$.
We we are going to prove the theorem in five cases.

$Case$ $1$. There exists a $x_i\in U$ such that $f(x_i)=y_1$ and $f(y_1)=x_j$.

We have $\delta_f(x_i)=t+n-2$ and $\delta_f(y_1)=t+n-2$.
Since $d(x_i,y_1)=d(f(x_i),f(y_1))=1$, $\delta_f(G_{(n,1)})\geq\delta_f(x_i)+\delta_f(y_1)=2t+2(n-2)$ and the equal sign holds if and only if $n=2$.
And in this case we do not need any restrictions other than condition $n=2$.
\\

$Case$ $2$. There exists a $x_i\in U$ such that $f(x_i)=y_1$ and $f(y_1)=y_j$.

It is easy to find that $\alpha(x_i,x_j)=\alpha(x_i,y_l)=(a_2,a_1)$ for $i\neq j$ and $l\neq1$, so we have $\delta_f(G_{(n,1)})\geq2t+2n-4$ by Lemma 2.2 and the equal sign holds if and only if $n=2$ and $f(U\setminus\{x_i\})\subseteq U$.
\\

$Case$ $3$. $f(x_i)\neq y_1$ for any $x_i\in U$ and $f(y_1)\in U$.

Without loss of generality, we can assume that $f(y_2)=y_1$.
Then $\alpha(y_1,x_i)=\alpha(y_1,y_j)=(a_1,a_2)$ for $j\neq2$, so $\delta_f(G)\geq 2(t+n-1)>2t$ by Lemma 2.2.
\\

$Case$ $4$. $f(x_i)\neq y_1$ for any $x_i\in U$ and $f(y_1)=y_1$.

Without loss of generality, we can assume that $f(y_i)\in U$ and $f(x_j)\in V(K_n)$ for $i\in\{2,3,\ldots,a+1\}$ and $j\in\{1,2,\ldots,a\}$.
Then $\delta_f(G)=2[(n-2)+(n-3)+\cdots+(n-a-1)]$.
Obviously, when $a=1$, $\delta_f(G)$ is going to be at least $2n-4$.
Let $2n-4\leq 2t$, then $n\leq t+2$.
\\

$Case$ $5$. $f(x_i)\neq y_1$ for any $x_i\in U$ and $f(y_1)=y_j$ for $j\neq1$.

Without loss of generality, we can assume that  $f(y_i)\in U$, $f(y_{a+2})=y_1$ and $f(x_j)\in V(K_n)$ for $i\in\{2,3,\ldots,a+1\}$ and $j\in\{1,2,\ldots,a\}$.
Then $\delta_f(y_{a+2})=t$ and $\delta_f(y_1)=t$.
Since $d(y_1,y_{a+2})=d(f(y_1), f(y_{a+2}))=1$, $\delta_f(G)\geq2t$ and the equal sign holds if and only if $a=0$ or $a=1$ and $n=3$.

From the above discussion, we know that
$$\pi(G)=
\left\{
    \begin{array}{lc}
        2n-4 & n \leq |U|+2=t+2 \\
        2t&otherwise\\
    \end{array}
\right.
$$

\end{proof}

Finally we talk about near automorphisms for $G_{(n,n)}$ and $n\geq 2$. \\
We also assume that $V(K_n)=\{y_1,y_2,\ldots,y_n\}$ and $U_i=\{x_{i1},x_{i2},\ldots,x_{it_i}\}$ is the suspension vertices set of $y_i$ for $i\in\{1,2,\ldots,n\}$.
We let $f(y_1)=y_2$, $f(y_2)=y_1$, $f(U_i)=U_i$ and $f(y_j)=y_j$ for $i\in\{1,2,\ldots,n\}$ and $j\in\{3,4,\ldots,n\}$, then we have $\sigma_f(G_{(n,n)})=2t_1+2t_2$.
And next we will prove $\pi(G)=t_1+t_2$.
\\
\begin{Theorem}
If $n\geq 2$ and $1\leq|U_1|=t_1\leq|U_2|=t_2\leq\cdots\leq|U_n|=t_n$, then $\pi(G_{(n,n)})=2t_1+2t_2$.
\end{Theorem}
\begin{proof}
$Case$ $1$. $n\geq 3$, $2\leq|U_1|=t_1\leq|U_2|=t_2\leq\cdots\leq|U_n|=t_n$ and $f(x_{ik})=y_j$ for some $x_{ik}\in U_i$.\\
Then $\alpha(x_{ik},x_{jl})=(a_3,a_2)$ or $(a_3,a_1)$ for $j\neq i$.
So we have $\delta_f(x_{ik},U_j)\geq t_j$ for $j\neq1$, and then $\delta_f(G)>2t_1+2t_2$ by $n\geq 3$ and Lemma 2.2.
So we have $f(U)=U$ for $U=U_1\cup U_2\cup\cdots\cup U_n$.
If $x_ij\in U_i$ and $f(x_{ij})\notin U_i$, by $Case$ $1$ of Lema 3.1, we have $U_i$ and $U_j$ such that $f(U_i)=M_{i_1}\cup M_{i_2}\cup\cdots \cup M_{i_h}$ where $M_{i_r}\subseteq U_{i_r}$ and $M_{j_r}\subseteq U_{j_r}$, and we also have (1), (2), (1') and (2') are right.
It is easy to know that there is a $x\in U_i$ such that $\alpha(x,y_i)=(a_1,a_2)$ by $h\geq2$, so $\delta_f(y_i,U_i)\geq 1$ and equality occurs if and only if $h=2$, $|M_{i_1}|=t_i-1$ and $f(y_i)=y_{i_1}$ or $h=2$, $|M_{i_2}|=t_i-1$ and $f(y_i)=y_{i_2}$.
If $t_i=h=3$ and $|M_{i_1}|=|M_{i_2}|=M_{i_3}=1$, then there are at least two vertices $x_{i1}$ and $x_{i2}$ of $U_i$ such that $\alpha(x_{i1},y_i)=\alpha(x_{i2},y_i)=(a_1,a_2)$.
And we have the same thing for $U_j$.
Combining the results of (1),(2),(1') and (2'), we have : \\
\begin{equation}
\delta_f(U_i)+\delta_f(y_i,U_i)\ge \begin{cases}
  t_i& \text{ if } h = 2, |M_{i_1}| = 1, |M_{i_2}| = t_i-1, f(y_i) = y_{i_2} \ or\\
 & \text{  } h = 2,|M_{i_2}| = 1,  |M_{i_1}| = t_i-1, f(y_i) = y_{i_1} ; (a)\\

  5& \text{ if } t_i = h = 3, |M_{i_1}| = |M_{i_2}| = |M_{i_3}| = 1 \ (b)
\end{cases}
\end{equation}
\\

\begin{equation}
\delta_f(U_i)+\delta_f(y_i,U_i)\ge \begin{cases}
  t_i& \text{ if } h = 2, |M_{j_1}| = 1, |M_{j_2}| = t_j-1, f(y_j) = y_{j_2} \ or\\
 & \text{  } h = 2,|M_{j_2}| = 1,  |M_{j_1}| = t_j-1, f(y_j) = y_{j_1} ; (a')\\

  5& \text{ if } t_i = h = 3, |M_{i_1}| = |M_{i_2}| = |M_{i_3}| = 1 \ (b')
\end{cases}
\end{equation}
\\

and
$
$$
\delta_f(G)\geq
\left\{
    \begin{array}{lc}
        2((a)+(a'))=2(t_i+t_j) \\
        2((a)+(b'))=2(t_i+5)>2t_1+6\\
        2((b)+(a'))=2(t_j+5)>2t_1+6\\
        2((b)+(b'))=20>6+6=2t_i+2t_j
    \end{array}
\right.
$
\\
If there is a $U_i$ such that $f(U_i)\neq U_i$, then $\delta_f(G)\geq 2t_1+2t_2$ and equality occurs if and only if $f(x_{ik_1})\in U_j$, $f(x_{jk_2})\in U_i$, $f(x_{ls})\in U_l$ and $f(y_l)=y_l$ for $x_{ls}\notin\{x_{ik_1},x_{jk_2}\}$ and $l\in\{1,2,\ldots,n\}$ or $f(U_i\setminus x_{ik_1})\in U_j$, $f(U_j\setminus x_{jk_2})\in U_i$, $f(y_i)=y_j$, $f(y_j)=y_i$, $f(x_{ls})\in U_l$ when $t_1=|U_i|\leq |U_j|=t_2$ and $f(y_r)=y_r$ for $x_{ls}\notin\{U_i\setminus x_{ik_1},U_j\setminus x_{jk_2}\}$ and $r\notin \{i,j\}$ when $t_1=|U_i|=|U_j|$

Now we let $f(U_i)=U_i$. Since $\delta_f(G)>0$, there must be two vertices where $y_i$ and $y_j$ satisfy $f(y_i)=y_k$ and $y_j=y_l$ for $i\neq k$ and $j\neq l$.
So $\delta_f(y_i)=t_i+t_k$ and $\delta_f(y_j)=t_j+t_l$.
Hence it is easy to know that $\delta_f(G)\geq 2t_1+2t_2$ and  equality occurs if and only if $f(y_i)=y_j$, $f(y_j)=y_i$ and $f(y_l)=y_l$ for $t_1=t_i\leq t_j=t_2$ and $l\notin\{i,j\}$.
\\

$Case$ $2$. $n\geq 3$, $1=|U_1|=|U_2|=\cdots=|U_k|$, $2\leq|U_{k+1}|\leq\cdots\leq|U_n|$ and $1\leq k< n$.\\
If there is some vertex $x_{i1}\in U_i$ such that $f(x_{i1})\in U_j$  for $1=|U_i|<|U_j|=t_j$.
Then there has to a vertex $x_{lk_1}\in U_l$, might as well set $x_{l1}$ where it is satisfied $f(x_{l1})=x_{i1}$ for $|U_l|\geq2$.
It is easy to prove that $\delta_f(G_{(n,n)})$ is minimized when $l=j$, $t_j=t_2$ and $f$ preserves the distance of all other vertex pairs $(x,y)$ for $(x,y)\notin\{(x_{i1},y_i),(x_{i1},y_j),(x_i1,x_j),(x_{j1},x_j),(x_{j1},y_i),(x_{j1},y_j):x_j\in U_j\setminus x_{l1}\}$ or $l=j$, $t_j=t_2$, $f(y_i)=y_j$, $f(y_j)=y_i$ and $f$ preserves the distance of all other vertex pairs $(x,y)$ for $(x,y)\notin\{(x_i1,x_j),(x_{j1},x_j),(x_j,y_i),(x_j,y_j):x_j\in U_j\setminus x_{l1}\}$, and the values of $\delta_f(G_{(n,n)})$ are $2t_2+2$ and $4t_2-4$, respectively.
Obviously, we only need to discuss the latter case.
We let $4t_2-4<2t_2+2$, then we have $t_2<3$, so $t_2$ can only be equal to $2$ and $4t_2-4=4$ by $2\leq t_2<3$, otherwise $4t_2-4\geq 2t_2+2$.
And the thing to notice is that $1=|U_1|$ and $2\leq|U_2|=t_2\leq\cdots\leq|U_n|=t_n$ whenever that minimum occurs. \\
Now we let $f(U_i)=U_i$ for $|U_i|=1$ and, $W=\{U_j:t_j\geq2\}$ and $W'=\{U_j:f(U_j)\neq U_j, t_j\geq2\}\neq\varnothing$.
If $|W|=1$, then there is a $y_i$ with $f(y_i)\in W$ by the assumption, so $\delta_f(G_{(n,n)})>2t_1+2t_2$. \\
If $|W|\geq3$, we know that $\delta_f(G_{(n,n)})>2t_1+2t_2$ by $Case$ $1$.
So we only need to discuss $W=W'=\{U_{n_1},U_n\}$ and it is easy to prove that $\delta_f(G_{(n,n)})>2t_2+2$, and equality occurs if and only if $n=3$, $t_2=t_3$ and there is exactly one pair of points where $(x_{2k_1},x_{3k_2})$ such that $\alpha(x_{2k_1},y_1)=(a_1,a_2)$ and $\alpha(x_{3k_1},y_2)=(a_1,a_2)$. 

Now we let $f(U)_i=U_i$ for $i\in\{1,2,\ldots,n\}$.
In this case, $f$ and $\delta_f(G_{(n,n)})$ is exactly the same as $Case$ $1$ when $f(U_i)=U_i$ for $i\in\{1,2,\ldots,n\}$.
\\

$Case$ $3$. $n\geq 3$, $1=|U_1|=|U_2|=\cdots=|U_n|$.\\
It is easy to prove that $\delta_f(G)\geq4$, equality occurs if and only if there are exactly two pairs of vertices $(x_{i1},x_{j2})$ and $(y_i,y_j)$ that satisfy $\alpha(x_{i1},y_i)=\alpha(x_{j1},y_j)=(a_1,a_2)$, $\alpha(x_{i1},y_j)=\alpha(x_{j1,y_i}=(a_2,a_1)$ and  $f$ preserves the distance of all other vertex pairs.
\\

$Case$ $4$. $n=2$.\\
Since $G_{(1,1)}=P_4$ when $t_1=t_2=1$, we let $t_2\geq2$ and $U=U_1\cup U_2$.
Firstly, we prove that $\delta_f(G_{(n,n)})>2t_1+2t_2$ if $f(x_1)=y_1$, $f(x_2)=y_2$ for $x_1 ,x_2\in U$.
We are only going to prove it for the case where $x_1\in U_1$ and $x_2\in U_2$, and we can do the same for the other cases.
Since $\alpha(x_1,x_{2i})=(a_3,a_2)$ or $(a_3,a_1)$ and $\alpha(x_2,x_{1j})=(a_3,a_1)$ or $(a_3,a_1)$, $\delta_f(G_{(1,1)})>2t_1+2t_2$ by Lemma 2.1 and Lemma 2.2.
If $f(U)\neq U$, there is exactly one vertex $x\in U$ such that $f(x)=y_1$ or $y_2$. \\
If $x\in U_1$ $f(x)=y_1$ and $f(y_1)\in U_1$, then $f(y_2)=y_2$ from the above discussion, and it is easy to observe that $\delta_f(G_{(n,n)}\geq 2t_1+2t_2$ and equality occurs if and only if $f(U_1\setminus x)\subseteq U_1$ and $f(U_2)=U_2$. \\
If $x\in U_1$ $f(x)=y_1$ and $f(y_1)\in U_2$, it is easy to prove $\delta_f(G_{(1,1)})>2t_1+2t_2$. \\
We also have the same classification for the other cases when $f(U)\neq U$, and we have $\delta_f(G_{(1,1)})\geq2t_1+2t_2$ in every case as well. \\
Next we let $f(U)=U$ and $f(U_i)\neq U_i$ for $i\in\{1,2\}$.
Then it is easy to know that $\delta_f(G_{(1,1)})\geq2t_1+2t_2$ and equality occurs if and only if there is exactly one pair of vertices $(x_{1i},x_{2j})$ such that $\alpha(x_{1i},y_1)=(x_{2j},y_2)=(a_1,a_2)$. \\
For $f(U_i)=U_i$, $f$ can only be equal to $(y_1y_2)$, so $\delta_f(G)=2t_1+2t_2$.
To sum up, $\delta_f(G_{(1,1)})=2t_1+2t_2$.

\end{proof}


\begin{thebibliography}{s2}

\bibitem{AW}  W. Aitken, Total relative displacement of permutations, J. Combin. Theory (A) (1999) 1-21.
\bibitem{CF} C.-F. Chang, H.-L. Fu, Near automorphisms of trees with small total relative displacements, J.
Comb. Optim. 14 (2007) 191-195.
\bibitem{CCF} C.-F. Chang, B.-L. Chen, H.-L. Fu, Near automorphisms of cycles, Discrete Mathematics 308 (2008)
1088-1092.
\bibitem{HV} G. Chartrand, H. Gavlas, D.W. Vander Jagt, Near-automorphisms of graphs, Graph Theory, Combinatorics and Applications, in:Y. Alavi, D. Lick, A.J. Schwenk (Eds.), Proceedings of the 1996
Eighth Quadrennial International Conference on Graph Theory, Combinatorics, Algorithms and
Applications, vol. 1, New Issues Press, Kalamazoo, 1999, pp. 181-192.
\bibitem{R} K.B. Reid, Total relative displacement of vertex permutations of Kn1,n2,...,nt
, J. Graph Theory 41
(2002) 85-100.

\end{thebibliography}
\end{document}